\documentclass{article}

\usepackage{amsmath,amssymb}
\usepackage{hyperref}

\def\cA{{\cal A}}

\def\RP{{\mathbf{R}^+}}

\begin{document}



\title{On weighted estimates for a class of Volterra integral operators\footnote{Originally published in 
Doklady Ros. Akad. Nauk (1997), vol.~357, p.~455.
Translated in Doklady Mathematics (1997), vol.~56, no. 3 , p.~906-908. \url{https://www.elibrary.ru/item.asp?id=13256293}}}
\author{Vyacheslav S. Rychkov\\\it Moscow Institute of Physics and Technology}

\date{}
\maketitle

\begin{abstract}
Volterra integral operators ${\cal A}=\sum_{k=0}^m {\cal A}_k$, $({\cal A}_k f)(x)= a_k (x)\int_0^x t^k f(t) \,dt$, are studied acting between weighted $L_2$ spaces on $(0,+\infty)$. Under certain conditions on the weights and functions $a_k$, it is shown that $\cal A$ is bounded if and only if each ${\cal A}_k$ is bounded. This result is then applied to describe spaces of pointwise multipliers in weighted Sobolev spaces on $(0,+\infty)$.
\end{abstract}
{\bf 1.} In the theory of functions, it is rather common to study the possibility of weighted estimates of the kind
$$
 \| v {\cal A}f \|_p \le c \|  uf \|_p  \eqno  (1)
$$
for Volterra integral operators $\cA$, defined by the formula
$$
  ({\cal A}f)(x)=\int_0^x A(x,t)f(t) \,dt.
$$

Here the function $f$ and its image $\cA f$ are defined on the half-line $\RP =(0,+\infty)$; $u,v$ are nonnegative on $\RP$ functions (weights); $\|\cdot\|_p=\|\cdot\|_{L_p(\RP)}, 1<p<+\infty$; constant~$c>0$ is independent of $f$.

The answer (whether (1) takes place or not) must be given in terms of the kernel $A(x,t)$ and the weights $u,v$. 
At present this answer is known only for certain classes of non-negative kernels close in some sense to the kernel of the Riemann-Liouville operator $(x-t)^\alpha,\alpha\ge0$, and imposing at the same time minimal restrictions on the weights $u,v$. It seems that the most general class of such kernels was pointed out by R. Oinarov \cite{Oin93}.

This work considers a new class of kernels of the kind
$$
  A(x,t)=\sum_{k=0}^m a_k (x)t^k , \quad m \in \text{\bf N}.
  \eqno (2)
$$
Restrictions on functions $a_k$ are minimal or absent. In particular, these functions are not assumed related and can take values of either sign. Therefore, our class contains sign-changing kernels. On the other hand, estimate (1) for operators with kernels of the kind (2) is studied only for $p=2$. Some rather strong condition is also imposed on the weight $u$.

{\bf 2.} Denote by $B_\delta,\delta\ge0$, the set of positive locally integrable on $\RP$ functions~$w$ satisfying with some constant~$C_w$ the integral doubling condition
$$
\int_\Delta w(x) \,dx \le C_w \int_{\frac{1}{2}\Delta} w(x)\,dx
$$
for any interval $\Delta\subset\RP$ of length $|\Delta|\ge\delta$,
where $\frac{1}{2}\Delta$ is a twice smaller interval centered at the same point.

The main result of this work is the following theorem.

{\bf Theorem 1.}
{\it Let $u^{-2}\in B_\delta$ for some $\delta\ge 0$.
If $\delta>0$, then assume in addition $a_k v \in L_2 (0,r)\,\forall r>0$,$k=0\ldots m-1$.
Then to have the estimate
$$
 \| v {\cal A}f \|_2 \le c \|  uf \|_2,  \eqno  (3)
$$
where $\cA$ is an operator with the kernel (2), it is necessary and sufficient that
$$
s_k=\sup_{r>0} \|a_k v\|_{L_2 (r,+\infty)}\cdot
		   \|x^k u^{-1}\|_{L_2 (0,r)} <+\infty,\quad  k=0\ldots m.
		   \eqno  (4)
$$
}

{\bf Remarks.} 1) Denote by $L_{2,u}$ the weighted space of functions $f$ on $\RP$ with norm $\|fu\|_2$. Having (3) now means that ${\cal A}:L_{2,u}\to L_{2,v}$. Represent~$\cal A$ as a sum
${\cal A}=\sum_{k=0}^m {\cal A}_k$, where
$$
({\cal A}_k f)(x)=a_k(x)\int_0^x t^kf(t)\,dt.
$$
Inequality (3) for the operator ${\cal A}_k$ instead of $\cA$ reduces, via the substitution
$f_1(x)=x^kf(x),
u_1(x)=x^{-k}u(x), v_1(x)=a_k(x)v(x)$, to the well-studied weighted Hardy inequality
$$
\left\|v_1\int_0^xf_1(t)\,dt\right\|_2\le c\|u_1f_1\|_2,
$$
criterion for whose validity is known since a long time (see \cite{Maz'ya}) and takes the form
$$
\sup_{r>0}\|v_1\|_{L_2(r,+\infty)}\cdot\|u_1^{-1}\|
_{L_2(0,r)}<+\infty.
$$
From here it follows that
$$
s_k<+\infty\Longleftrightarrow{\cal A}_k:L_{2,u}\to L_{2,v}.
$$
Therefore, Theorem 1 is equivalent to the statement that (under its assumptions) we have ``splitting'' for the operator ${\cal A}=\sum_{k=0}^m {\cal A}_k$, in the sense that
$$
{\cal A}:L_{2,u}\to L_{2,v}\Longleftrightarrow{\cal A}_k:L_{2,u}\to L_{2,v},
\quad k=0\ldots m.
$$

2) Our method of proof of Theorem 1 gives the following estimate for the smallest constant $c$ in inequality (3) (or, what is the same, the norm of the opertator $\cA$ acting from $L_{2,u}$ to $L_{2,v}$):
$$
c_1\sum_{k=0}^ms_k\le\|\cA\|_{L_{2,u}\rightarrow L_{2,v}}\le
c_2\sum_{k=0}^ms_k.
$$
Constant $c_2$ here is universal. As for $c_1$, this constant
depends on $m$, on $C_{u^{-2}}$ (the doubling constant for $u^{-2}$), as well as
(if $\delta>0$) on the quantity
$\sum_{k=0}^{m-1} \|a_k v\|_{L_2 (0,r_0)}\cdot\|x^ku^{-1}\|_{L_2 (0,r_0)}$,
where $r_0$ is determined by $m,C_{u^{-2}},\delta$.

{\bf 3.} The following statement plays the main role in the proof of Theorem 1.

{\bf Lemma 1.}{\it Let $u^{-2}\in B_\delta,m\in{\hbox{\bf N}}$. Then there exists
such $r_0\ge0$ ($r_0=0$ for $\delta=0$), that
$$
\hbox{{\bf G}}(u^{-1}\chi_r,xu^{-1}\chi_r,\ldots,x^mu^{-1}\chi_r)\ge
$$
$$
\ge\varepsilon\|u^{-1}\chi_r\|_2\|xu^{-1}\chi_r\|_2\cdots
\|x^mu^{-1}\chi_r\|_2,
$$
$$
r\ge r_0,
$$
with a constant $\varepsilon>0$ independent of $r\ge r_0$.
}

Here $\hbox{\bf G}$ is the Gram determinant of a system of functions in $L_2$,
$\chi_r=\chi_{(0,r)}$ is a characteristic function of the interval.
This statement therefore asserts the uniform in $r$ non-degeneration of the parallelepiped with edges
$u^{-1}\chi_r,\ldots,x^mu^{-1}\chi_r$.

Let us show how to derive Theorem 1 from this Lemma. From the above discussion it's clear that we only need to prove the necessity of conditions (4). Inequality (3) written for a function $f_r$ supported on $[0,r]$ implies the inequality
$$
\left \| v(x)a_0(x)\int_0^rf_r(t)dt+
  v(x)\sum_{k=1}^ma_k(x)\int_0^rt^kf_r(t)\,dt \right\|_{L_2(r,+\infty)}
  \le c \|  uf_r \|_{L_2(0,r)}. \eqno (5)
$$
Suppose $uf_r$ belongs to the orthogonal complement of the linear span
$E_r$ of the set of functions $xu^{-1}\chi_r,\ldots,x^mu^{-1}\chi_r$.
By Lemma 1, the angle between the vector $u^{-1}\chi_r$ and the subspace
$E_r$ is separated from zero uniformly in
$r\ge r_0$. Therefore we can choose $f_r$ so that the angle between
$uf_r$ and $u^{-1}\chi_r$ is uniformly separated from $\pi/2$, i.e.
$$
\int_0^r f_r(t)\,dt\ge\alpha
\|uf_r\|_{L_2(0,r)}\cdot\|u^{-1}\|_{L_2(0,r)},\quad r\ge r_0\quad (\alpha>0).
$$
In this case (5) implies
$$
\|a_0v\|_{L_2(r,+\infty)}\cdot\|u^{-1}\|_{L_2(0,r)}\le c/\alpha,
\quad r\ge r_0,
$$
from where $s_0<+\infty$. Finiteness of the other constants $s_k$ is shown by induction.

We see that the argument uses essentially the geometry of the Hilbert space $L_2$.

{\bf 4.} Let us consider some applications. Consider on $\RP$ the weighted Sobolev space
$W=W_{2,u}^{(l)}$ with the norm $\|f\|_W=\|f\|_{L_2(0,1)}+
\|f^{(l)}u\|_2$. Particularity of this norm is that by taking the norm of the function itself only on an initial interval of $\RP$ allows to include into this space polynomials of degree $\le l-1$. Spaces $W_{2,u}^{(l)}$
were introduced and studied by L.D. Kudryavtsev in \cite{Kud}, using an equivalent norm
$\sum_{k=0}^{l-1}|f^{(k)}(0)|+\|f^{(l)}u\|_2$.

Function $\varphi$ is called (pointwise) multiplier from ${}^1W$
into ${}^2W$, if $\varphi f\in {}^2W\;\forall f\in{}^1W$. The space of multipliers is denoted by $M({}^1W\to {}^2W)$.

Various aspects of the theory of multipliers in unweighted spaces of differentiable functions were studied in the book \cite{MS}. G.A.Kalyabin \cite{Kalyabin} described multipliers in Sobolev spaces on $\hbox{{\bf R}}^n$ 
with the norm
$\|f\|_{L_p(B(0,1))}+\|\nabla_l f\|_p$ in the case $p>n$; he also posed the question of describing
multipliers in the considered here weighted case.

Obtained in Theorem 1 criterion of boundedness in weighted spaces of operators with kernels of the kind (2) is decisive to prove the following result.

{\bf Theorem 2.}{\it Let $u^{-2}\in B_\delta$, $v^{-1}\in L_2(0,r)
\,\forall r>0$. Then the space $M(W_{2,u}^{(l)}\to W_{2,v}^{(m)}),
m\le l$, consists of those and only those $\varphi$ satisfying the following two conditions:
$$
\|(\varphi x^k)^{(m)}v\|_2<+\infty,\quad k=0\ldots l-1 \eqno(6)
$$
$$
\sup_{r>0}\|(\varphi x^k)^{(m)}v\|_{L_2(r,+\infty)}\cdot
\|x^{l-k-1}u^{-1}\|_{L_2(0,r)}<+\infty,
\quad k=0\ldots l-1, \eqno(7)
$$
to which in the case $m=l$ one more condition is added, namely
$$
\|\varphi v u^{-1}\|_{L_\infty(\RP)}<+\infty. \eqno(8)
$$

}
{\bf Remark.} For
$m=l,u=v,(1+x^{l-1})u^{-1}\in L_2(\RP)$
the considered space of multipliers was described by the author in \cite{Rychkov}. The obtained result had the form of the combination of two conditions
(6),(8), which is natural, since in that case (6)$\Rightarrow$(7).
Theorem 3 allows to widen significantly the class of weights for which a description of multipliers is available. Interestingly, multipliers in spaces with exponentially decreasing at $\infty$ weights remain not studied: their description should become the subject of future investigations.

Two formulas from the following lemma show how operators with kernels of the form (2) appear in the problem about multipliers.

{\bf Lemma 2.} {\it Let function $g$ on $\RP$ be such that
$$
g^{(k)}(0)=0,\quad k=0\ldots l-1.
$$
Then
\begin{equation}
(\varphi g)^{(l)}(x)= \varphi(x)g^{(l)}(x)+\frac1{(l-1)!}\sum_{k=0}^{l-1}C_{l-1}^k \left(\varphi x^k\right)^{(l)}\int_0^x(-t)^{l-k-1}g^{(l)}(t)\,dt;
\end{equation}
\begin{equation}
(\varphi g)^{(m)}(x)= \frac1{(l-1)!}\sum_{k=0}^{l-1}C_{l-1}^k \left(\varphi x^k\right)^{(m)}\int_0^x(-t)^{l-k-1}g^{(l)}(t)\,dt, \quad m<l.
\end{equation}
}

Detailed proofs of all statements will be published in Proceedings of Steklov Institute of Mathematics.

{\bf Note added (September 1997):} In the time since submission of this article, proofs of the given results have appeared in \cite{Rychkov2}. The authors has also obtained generalizations for $p\ne2$ \cite{Rychkov3}.


\begin{thebibliography}{9}
\bibitem{Oin93} R. Oinarov, “Two-sided norm estimates for certain classes of integral operators”, in {\it Investigations in the theory of differentiable functions of many variables and its applications. Part 16}, Trudy Mat. Inst. Steklov., 204, Nauka, Moscow, 1993, 240–250; [Translated in: Proc. Steklov Inst. Math., 204 (1994), 205–214] \url{http://mi.mathnet.ru/eng/tm1271}

\bibitem{Maz'ya} V.G. Maz'ya, {\it Sobolev Spaces}, Springer, 1985.

\bibitem{Kud} L.D.~Kudryavtsev, ``On norms in weighted spaces of functions given on infinite intervals',' Analysis Mathematica 12, 269–282 (1986) \url{https://doi.org/10.1007/BF01909365}

\bibitem{MS} V. G. Maz'ya and
T. O. Shaposhnikova, {\it Theory of multipliers in spaces of differentiable functions}, Monographs and Studies in Mathematics, vol. 23, Pitman Publishing Co., Brooklyn, New York, 1985

\bibitem{Kalyabin} G.A. Kalyabin, “Pointwise multipliers in some Sobolev spaces containing unbounded functions”, in {\it Investigations in the theory of differentiable functions of many variables and its applications. Part 16}, Trudy Mat. Inst. Steklov., 204, Nauka, Moscow, 1993, 160–165; [Translated in: Proc. Steklov Inst. Math., 204 (1994), 137–141] \url{http://mi.mathnet.ru/tm1266}

\bibitem{Rychkov} V. S. Rychkov, “Pointwise multiplicators in weighted Sobolev spaces on a half-line”, Math. Notes, 56:1 (1994), 704–710 \url{https://doi.org/10.1007/BF02110561}

\bibitem{Rychkov2} V. S.~Rychkov, “Splitting of Volterra integral operators with degenerate kernels”, in {\it Investigations in the theory of differentiable functions of many variables and its applications. Part 17, Collection of articles}, Tr. Mat. Inst. Steklova, 214, Nauka, Moscow, 1997, 267–285; [Translated in: Proc. Steklov Inst. Math., 214 (1996), 260–278] \url{http://mi.mathnet.ru/eng/tm1040}

\bibitem{Rychkov3} V.S.~Rychkov, ``Some Weighted Hardy-Type Inequalities and Applications,'' Proc. of A. Razmadze Georgian Math. Inst. (1997), vol.~112. p.~113-129 \url{http://www.rmi.ge/proceedings/volumes/ps/v112-5.ps.gz}.

\end{thebibliography}
\end{document}